\documentclass[11pt]{article}
\usepackage{graphicx,amssymb,mathrsfs,amsmath,color}
\newfont{\bb}{msbm10}

\newcommand{\tr}{^{\sf T}}

\def\Diag{{\rm Diag}}

\newcommand{\m}[1]{{\bf{#1}}}

\newcommand{\C}[1]{{\cal {#1}}}

\def\be{\begin{eqnarray}}
\def\ee{\end{eqnarray}}
\def\ben{\begin{eqnarray*}}
\def\een{\end{eqnarray*}}
\def\ba{\begin{array}}
\def\ea{\end{array}}

\newtheorem{theorem}{Theorem}[section]

\newtheorem{lemma}{Lemma}[section]
\newtheorem{corollary}{Corollary}[section]

\begin{document}
\cleardoublepage \pagestyle{plain}
\bibliographystyle{plain}

\title{ Convergence Revisit on Generalized Symmetric ADMM
\thanks{The work was supported by the National Natural  Science Foundation of China (Nos. 11671318; 11571271; 11631013) and
the Natural  Science Foundation of Fujian Province (No. 2016J01028). The second author  Xiaokai Chang was supported by the Hongliu Foundation of First-class Disciplines of Lanzhou University of Technology.}
   }

\author{
Jianchao Bai
\footnote{Department of Applied Mathematics,   Northwestern Polytechnical
University, Xi'an, 710129,   China.  Past addresses: School of Mathematics and Statistics, Xi'an Jiaotong University, Xi'an 710049,   China
(\tt bjc1987@163.com).}
\quad
Xiaokai Chang
\footnote{College of Science, Lanzhou University of Technology, Lanzhou  730050,   China
(\tt xkchang@lut.cn).}
\quad Jicheng Li
\footnote{School of Mathematics and Statistics, Xi'an Jiaotong University, Xi'an 710049,   China
(\tt jcli@mail.xjtu.edu.cn).}
\quad Fengmin Xu
\footnote{School of Economics and Finance, Xi'an Jiaotong University, Xi'an 710049,   China
(\tt fengminxu@mail.xjtu.edu.cn).}
}
\date{}
\maketitle

\centerline{\small\it\bf Abstract}\vskip 1mm
In this note, we show a sublinear nonergodic convergence rate for the algorithm developed in [Bai, et al. Generalized symmetric ADMM for separable convex optimization. Comput. Optim. Appl. 70, 129-170 (2018)], as well as its linear convergence   under   assumptions that the sub-differential of each component objective function   is piecewise linear and all the constraint sets are polyhedra. { These remaining convergence results  are established for the stepsize  parameters of dual variables belonging to a special isosceles triangle region, which aims to  strengthen  our understanding for convergence   of the generalized symmetric ADMM.}


\bigskip

\section{Introduction}
Revisit the following prototype multi-block separable convex optimization
\begin{equation} \label{Prob}
\begin{array}{lll}
\min  &  \sum\limits_{i=1}^{p}f_i(x_i)+ \sum\limits_{j=1}^{q}g_j(y_j) \\
\textrm{s.t. } &   \sum\limits_{i=1}^{p}A_i x_i+\sum\limits_{j=1}^{q}B_jy_j=c,\\
     &   x_i\in \mathcal{X}_i, \; i=1, \cdots, p, \\
     &   y_j\in \mathcal{Y}_j, \; j =1, \cdots, q, \\
\end{array}
\end{equation}
 where $f_i(x_i):\mathbb{R}^{m_i}\rightarrow{ \mathbb{R}},  g_j(y_j):\mathbb{R}^{d_j}\rightarrow\mathbb{R}$ are closed and proper convex functions (possibly nonsmooth); $A_i\in\mathbb{R}^{n\times m_i}, B_j\in\mathbb{R}^{n\times d_j}$ and $ c\in\mathbb{R}^{n}$ are given  matrices and vectors, respectively;   $\mathcal{X}_i\subset \mathbb{R}^{m_i}$ and $\mathcal{Y}_j\subset \mathbb{R}^{d_j} $   are  polyhedra; $p \ge 1$ and $q \ge 1$ denote two integers. Throughout we assume   the solution set of the problem (\ref{Prob}) is nonempty and all the matrices $A_i(i=1,\cdots, p)$ and $B_j( j=1,\cdots,q)$ have full column rank.

 By denoting { $\mathcal{A}=\left[A_1, \cdots,A_p\right],\mathcal{B}=\left[B_1, \cdots,B_q\right],
\m{x} = (x_1, \cdots, x_p)$ and $\m{y} = (y_1, \cdots, y_q) $,}
  the augmented Lagrangian function of the problem (\ref{Prob}) is written as
\[
 \mathcal{L}_\beta\left(\m{x}, \m{y},\lambda\right)=L\left(\m{x}, \m{y},\lambda\right)+\frac{\beta}{2} \left\| \C{A} \m{x} + \C{B} \m{y} -c\right\|^2,
\]
where $\beta>0$ is a penalty parameter   and
\[
L\left(\m{x}, \m{y},\lambda\right)=\sum\limits_{i=1}^{p}f_i(x_i)+ \sum\limits_{j=1}^{q}g_j(y_j) -  \langle \lambda, \C{A} \m{x} + \C{B} \m{y} -c \rangle
\]
denotes the Lagrangian function  associated with a Lagrange multiplier $\lambda\in \mathbb{R}^{n}$.
As studied  in our recent work \cite{BaiLiXuZhang2017}, the  Generalized Symmetric Alternating Direction Method of Multipliers  (\textbf{GS-ADMM})    reads the following updates:
\begin{equation}\label{GS-ADMM}
   \left\{\begin{array}{lll}
 \textrm{For}\ i=1,2,\cdots,p,\\
 \quad x_{i}^{k+1}=\arg\min\limits_{x_{i}\in\mathcal{X}_{i}} \mathcal{L}_\beta (x_1^k,\cdots, x_i,\cdots,x_p^k,\m{y}^k,\lambda^k)
+P_i^k(x_i), \\
\quad \textrm{where } P_i^k(x_i) =
\frac{\sigma_1\beta}{2}\left\|A_{i}(x_{i}-x_{i}^k)\right\|^2,\\
 \lambda^{k+\frac{1}{2}}=\lambda^k-\tau\beta(\mathcal{A} \m{x}^{k+1}+\mathcal{B}\m{y}^{k}-c),\\ \\
 \textrm{For}\ j=1,2,\cdots,q,\\
 \quad y_{j}^{k+1}=\arg\min\limits_{y_{j}\in\mathcal{Y}_{j}} \mathcal{L}_\beta(\m{x}^{k+1},y_1^k,\cdots, y_j,\cdots,y_q^k,\lambda^{k+\frac{1}{2}}) + Q_j^k(y_j), \\
 \quad \textrm{where } Q_j^k(y_j) =
\frac{\sigma_2\beta}{2}\left\|B_{j}(y_{j}-y_{j}^k)\right\|^2,\\
\lambda^{k+1}=\lambda^{k+\frac{1}{2}}-s\beta(\mathcal{A}\m{x}^{k+1}+\mathcal{B}\m{y}^{k+1}-c),
\end{array}\right.
\end{equation}
where $\tau$ and $s$ are   stepsize parameters  satisfying
\[
  (\tau, s) \in \C{G} = \left\{ (\tau, s) \ | \ \tau + s >0,\ -\tau^2 - s^2 -\tau s + \tau + s + 1 >0 \right\},
\]
and $\sigma_1\in (p-1,+\infty),\sigma_2\in (q-1,+\infty)$ are proximal parameters
for the regularization terms $P_i^k(\cdot)$ and $Q_j^k(\cdot)$, respectively.
\begin{figure}[htbp]\label{figG}
 \begin{minipage}{1\textwidth}
 \centering
\resizebox{14cm}{6.2cm}{\includegraphics{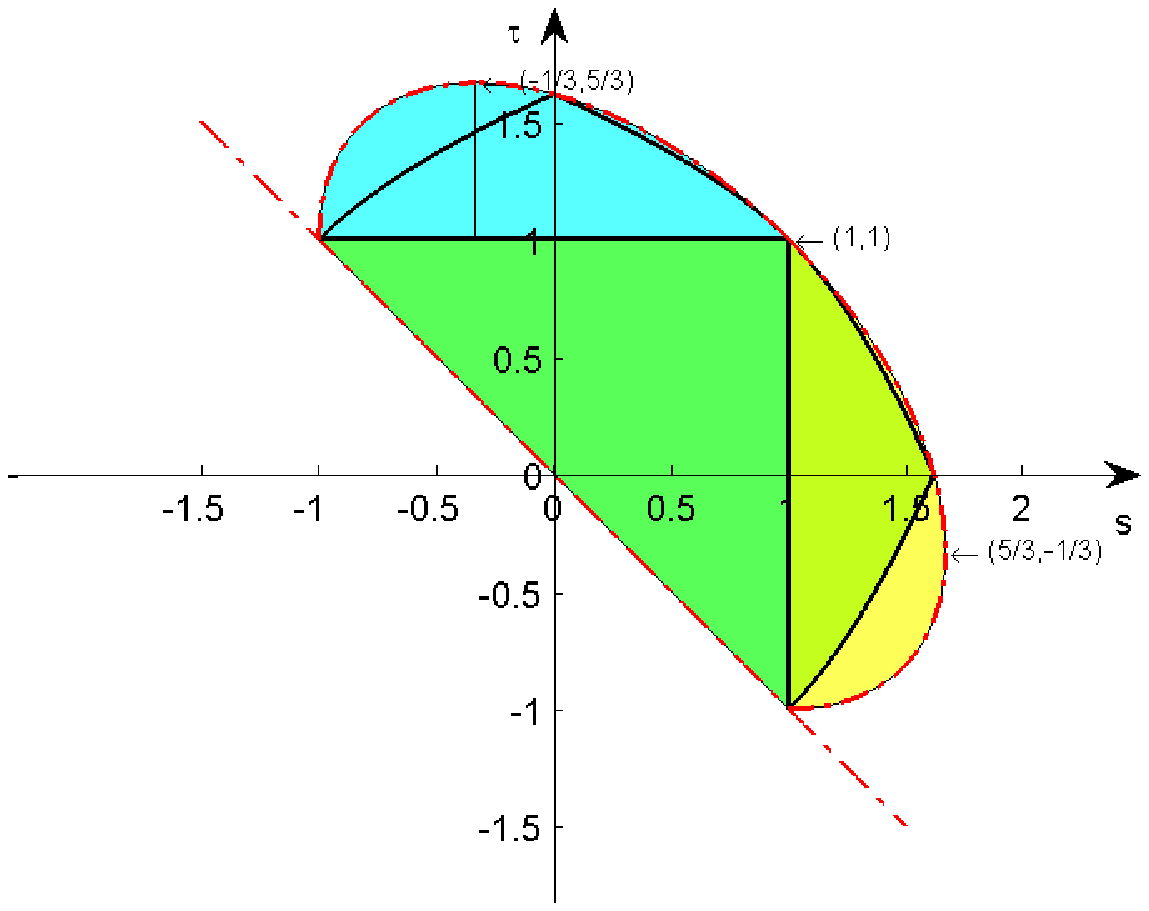}\includegraphics{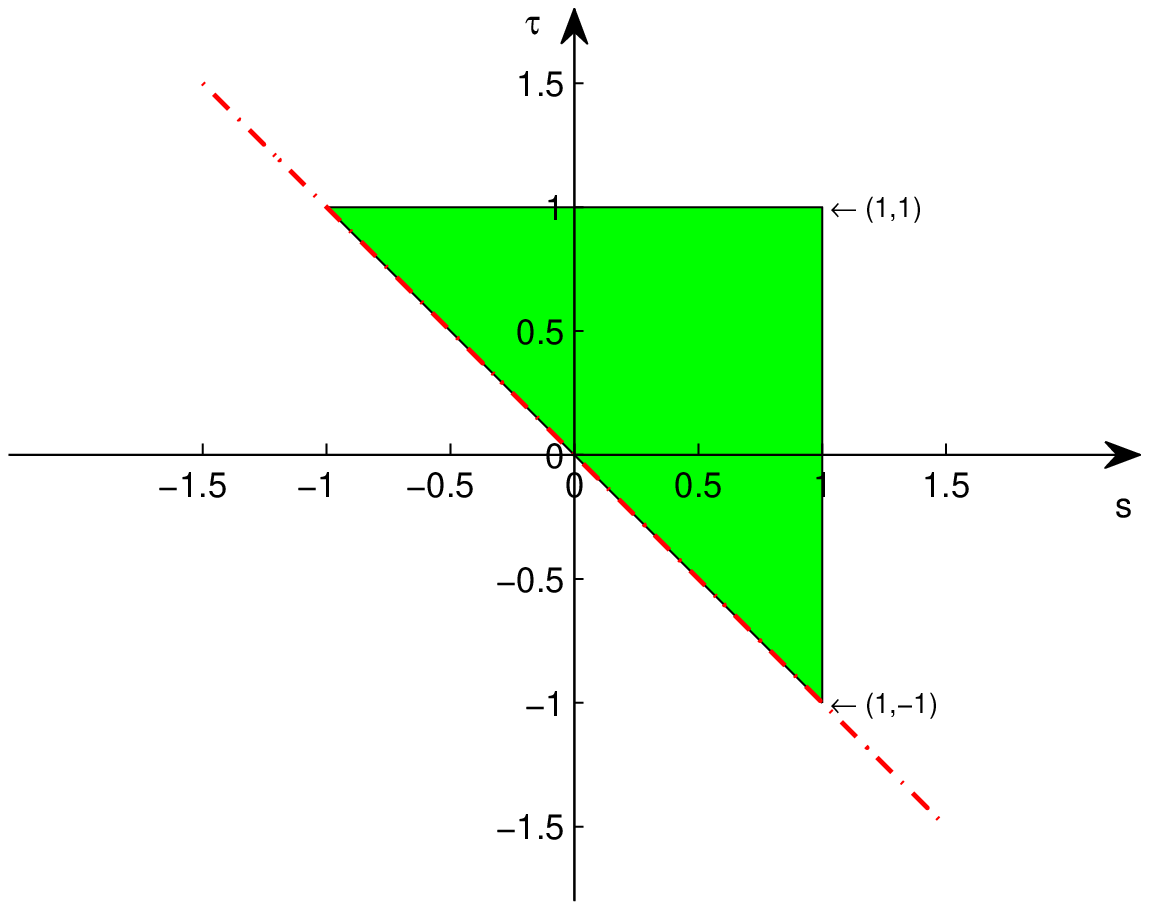}}
   \end{minipage}
\begin{center}
{\footnotesize Fig. 1: The left is the region $\mathcal{G}$ and the right is the region $\mathcal{D}$}
\end{center}
\end{figure}
By making use of  a prediction-correction interpretation for GS-ADMM, we analyzed its global convergence, sublinear convergence rate in the ergodic sense and   convergence complexity of two special cases allowing either $\sigma_1$ or $\sigma_2$ to be zero. However,  two  remaining tasks   were not  settled as mentioned by the past reviewers: (1)
  How to establish its worst-case $\mathcal{O}(1/t)$ convergence rate in the nonergodic sense, where $t$ denotes the iteration number?
 (2)
   Whether there exists a linear convergence rate of GS-ADMM under  some  mild  assumptions?
This note aims to  give   positive answers for these   questions but for the following subregion (shown in the right-hand side of Fig. 1) of $\mathcal{G}$, that is,
\begin{equation}\label{setK}
(\tau, s) \in \mathcal{D}:= \left\{(\tau, s) \ | \ \tau<1,\ s<1,\ \tau+s>0 \right\}.
\end{equation}
Notice that the above region   is   much wider than that
($\tau=s\in(0,1)$) in   \cite[Algorithm 3]{Heyuan2015}. Moreover,  it can be seen   by later analysis  that  the symmetric ADMM
(S-ADMM, \cite{HeMaYuan2016}) for solving the two-block separable convex optimization
also has the worst-case $\mathcal{O}(1/t)$  convergence rate
in the nonegodic sense as well as global linear convergence
rate { for parameters belonging to $\mathcal{D}$.}

{
\subsection{Relationship of GS-ADMM  to  related works}
The algorithm GS-ADMM was initially  proposed to generalize the meaningful S-ADMM \cite{HeMaYuan2016}  for solving the grouped multi-block separable convex optimization  problem (\ref{Prob}), whose convergence and iteration complexity could be still ensured for a larger domain of stepsizes  of dual variables than that introduced in \cite{HeMaYuan2016}. In practise, convergence   of GS-ADMM was analyzed by estimating the lower bound of $\left\|\m{w}^{k}-\widetilde{\m{w}}^k\right\|_G^2$ directly and by treating the domain of stepsize parameters as a whole,  while convergence of S-ADMM was showed separately by splitting the   domain of $(\tau, s)$ into several subdomains, where   $\m{w}^{k+1}$ and $\widetilde{\m{w}}^k$ are    called the predictive
variable  and the correcting variable, respectively. Note that  by taking $\sigma_1=\sigma_2=0$,   GS-ADMM with $p=q=1$ will become S-ADMM but continue to converge in the relatively larger convergence domain $\C{G}$. In addition, the original S-ADMM only works for
the two-block case and may not be convenient for solving large-scale   problems, while GS-ADMM could handle   large-scale multiple block problems since the block variables within each group
were updated in a Jacobian scheme.

Regardless of the additional dual variable update $\lambda^{k+\frac{1}{2}}$ (i.e. $\tau=0$), then GS-ADMM becomes a proximal ADMM-type algorithm with $s\in(0,\frac{1+\sqrt{5}}{2})$. Moreover, it will become the classical ADMM proposed by Glowinski-Marrocco \cite{GlowinskiMarrocco1975} when considering the simple two block case without using proximal regularization terms. To the best of knowledge, the first proximal ADMM was proposed by Eckstein \cite{Eckstein94} as GS-ADMM with $p=q=1, (\tau,s)=(0, 1)$  and with the following proximal terms
\[
P_1^k(x_1) =
\frac{1}{2}\left\|x_{1}-x_{1}^k\right\|_{\mathcal{T}_1}^2, \quad
Q_1^k(y_1) =
\frac{1}{2}\left\|y_{1}-y_{1}^k\right\|_{\mathcal{T}_2}^2,
\]
where  $\mathcal{T}_i=\frac{\mu_i^2}{\beta}I $ for any nonzero  scalars $\mu_i, i=1,2.$  Later,
a perfect extension on convergence analysis from the classical ADMM to  GS-ADMM   with $p=q=1$ and $ \tau=0$, but allowing the stepsize $s$ to stay in the   range $(0,\frac{1+\sqrt{5}}{2})$ was studied, see   Xu-Wu \cite{XuWu11} and Fazel, et. al. \cite{Fazelun13} for more details. Recently, He-Xu-Yuan \cite{HeY16} constructed a   proximal ADMM  for solving the problem (\ref{Prob}) with only $p$ block variables, and their algorithm could be regarded as a  special version of GS-ADMM with $(\tau,s)=(0,
1)$   barring the $y_j$-updates. Especially, the partially proximal ADMM-type algorithm \cite{SUn18} with a specified regularization term $Q_j^k(y_j)$ as ours   could be treated as the case that GS-ADMM with $p=1,\sigma_1=0$ and $\tau=1.$ Considering  the middle update $\lambda^{k+\frac{1}{2}}$(i.e. $\tau\neq 0$), convergence domain of the dual stepsizes of  GS-ADMM is still larger than that in the symmetric ADMM with indefinite proximal regularization \cite{gaoma18,SS18}.
}

\subsection{Notations  and  organizations}
Throughout the { note}, the symbols  $\mathbb{R}, \mathbb{R}^n, \mathbb{R}^{m\times n}$ denote the sets
of  real numbers,   $n$ dimensional real column vectors and   $m\times n$ real matrices, respectively. For any $x, y\in \mathbb{R}^n$,
$\langle x, y\rangle=x \tr y$ represents their inner product and
$\|x\|=\sqrt{\langle x, x\rangle}$ denotes the Euclidean norm of $x$, where
$\tr$ {  denotes the  transpose operation}. For any   symmetric matrix $G$, we define $\|x\|_G^2 = x \tr G x$ which is not necessarily nonnegative unless
$G$ is   positive definite.  The symbols $\lambda_{\max}(\cdot)$ and $\lambda_{\min}(\cdot)$   denote respectively the maximum  and minimum eigenvalue  of a square  matrix. The notations $I$ and $\bf{0}$     stand   for the identity matrix and   zero matrix with
proper dimensions, respectively. We call  $\phi(x)$   a piecewise linear multifunction if its graph $\{(x,y)|\ y\in \phi(x)\}$ is a union of finitely many polyhedra.
For  convenience,   let \[\mathcal{M}=\mathcal{X} \times \mathcal{Y} \times \mathbb{R}^n\] and the corresponding solution set be $\mathcal{M}^*,$ where $\C{X} = \C{X}_1 \times \C{X}_2 \times \cdots \C{X}_p $ and
$\C{Y} = \C{Y}_1 \times \C{Y}_2 \times \cdots \C{Y}_q $.
 We also preset
\[
\m{u}=\left(\begin{array}{c}
\m{x}\\  \m{y}\\
\end{array}\right),\
\m{w}=\left(\begin{array}{c}
\m{x}\\ \m{y}\\  \lambda
\end{array}\right), \
\mathcal{J}(\m{w})=\left(\begin{array}{c}
-\C{A}\tr\lambda\\
-\C{B}\tr\lambda\\
 \C{A} \m{x} + \C{B} \m{y} - c
\end{array}\right),
\]
\begin{equation}\label{tilde-xy}
\widetilde{\m{x}}^k=\left(\begin{array}{c}
\widetilde{x}_1^k\\  \widetilde{x}_2^k\\ \vdots\\ \widetilde{x}_p^k
\end{array}\right)=\left(\begin{array}{c}
x_1^{k+1}\\  x_2^{k+1}\\ \vdots\\ x_p^{k+1}
\end{array}\right),
\quad \widetilde{\m{y}}^k=\left(\begin{array}{c}
\widetilde{y}_1^k\\  \widetilde{y}_2^k\\ \vdots\\ \widetilde{y}_q^k
\end{array}\right)=\left(\begin{array}{c}
y_1^{k+1}\\  y_2^{k+1}\\ \vdots\\ y_q^{k+1}
\end{array}\right),
\end{equation}
\begin{equation}\label{tilde-uw}
\widetilde{\m{u}}^k=\left(\begin{array}{c}
\widetilde{\m{x}}^k\\  \widetilde{\m{y}}^k\\
\end{array}\right),\quad
\widetilde{\m{w}}^{k}=\left(\begin{array}{c}
 \widetilde{\m{x}}^k\\ \widetilde{\m{y}}^k\\ \widetilde{\lambda}^k
\end{array}\right)=
\left(\begin{array}{c}
\m{x}^{k+1}\\ \m{y}^{k+1}\\   \widetilde{\lambda}^{k}
\end{array}\right),
\end{equation}
and
\begin{equation} \label{tilde-lambda}
\widetilde{\lambda}^k=\lambda^k-\beta \left(\mathcal{A}\m{x}^{k+1}+\mathcal{B}\m{y}^{k}-c\right).
\end{equation}

The rest of this paper is organized as follows. In Section 2, by making use of some  well-known identities, inequalities and matrix decomposition techniques, we first establish  {  sublinear} convergence rate of GS-ADMM in the nonergodic sense. Then, its global linear convergence rate, measured by an error function  $\textrm{dist}^2_{H}(\m{w}^{k+1}, \mathcal{M}^*)$ or $\left\|\m{w}^k-\m{w}^\infty \right\|_{H}$, is analyzed under mild assumptions.   Finally, we briefly conclude the paper  in Section 3.

\section{Main results}
At the beginning of this section,  we  first  analyze the worst-case $\mathcal{O}(1/t)$ nonergodic convergence rate of GS-ADMM   for any $(\tau,s)\in \mathcal{D}$. Then,   by   using several well-known  inequalities its   convergence rate is strengthened to   linear under the assumption that the subdifferential of each objective function  is  {  piecewise linear}.

\subsection{Sublinear nonergodic  convergence rate}
Let us review the following two basic lemmas given in \cite{BaiLiXuZhang2017}, which aims to interpret the GS-ADMM into a prediction-correction procedure.

\begin{lemma}\label{tilde-VI}
For the iterates $\widetilde{\m{u}}^k, \widetilde{\m{w}}^k$ defined in (\ref{tilde-uw}),
we have $\widetilde{\m{w}}^k\in\mathcal{M}$ and
\begin{equation} \label{app-vi}
  h(\m{u})-h(\widetilde{\m{u}}^k)+ \left\langle \m{w}-\widetilde{\m{w}}^k, \mathcal{J}(\widetilde{\m{w}}^k)+Q(\widetilde{\m{w}}^k-\m{w}^k)\right\rangle \geq 0, \quad \forall \m{w}\in \mathcal{M},
\end{equation}
where
$h(\m{u})=\sum\limits_{i=1}^{p}f_i (x_i) +\sum\limits_{j=1}^{p}g_j (y_j)$ and
\begin{equation} \label{Q}
Q=\left[\begin{array}{cc}
H_{\mathbf{x}} & \bf{0}\\
\bf{0} & \widetilde{Q}
\end{array}\right]
\end{equation}
with
\begin{equation} \label{Hx}
H_{\mathbf{x}}=\beta\left[\begin{array}{ccccccc}
\sigma_1 A_1\tr A_1  && - A_1\tr A_2         && \cdots  &&- A_1\tr A_{p}  \\
- A_2\tr A_1         && \sigma_1 A_2\tr A_2  && \cdots  &&- A_2\tr A_{p}  \\
\vdots               && \vdots              &&\ddots   &&\vdots   \\
 - A_p\tr A_1        &&  - A_p\tr A_2        &&\cdots  &&\sigma_1 A_p\tr A_p
\end{array}\right],
\end{equation}
\[
 \widetilde{Q}=\left[\begin{array}{ccccccc|c}
 (\sigma_2+1)\beta B_{1}\tr B_{1}  && {\bf 0}  &&\cdots &&{\bf 0}  &-\tau B_{1}\tr \\
 {\bf 0} && (\sigma_2+1)\beta B_{2}\tr B_{2} &&\cdots &&{\bf 0} &-\tau B_{2}\tr \\
 \vdots && \vdots &&\ddots &&\vdots &\vdots \\
 {\bf 0} && {\bf 0} &&\cdots &&(\sigma_2+1)\beta B_{q}\tr B_{q} &-\tau B_{q}\tr \\ \hline
-B_{1} && -B_{2}  &&\cdots &&-B_{q} &\frac{1}{\beta}I
\end{array}\right].
\]
\end{lemma}

\begin{lemma} \label{pre-lem}
 For the sequences  $\{\m{w}^k\}$ and $\{\widetilde{\m{w}}^k\}$ generated by  GS-ADMM, the following equality  holds
\begin{equation} \label{pre-core}
\m{w}^{k+1}=\m{w}^k-M(\m{w}^k-\widetilde{\m{w}}^k),
\end{equation}
 where
\begin{equation} \label{pre-core-M}
M=\left[\begin{array}{c|cccc}
I &     &  & & \\ \hline
&   I &  & & \\
&    &\ddots  & & \\
&    &  &I & \\
&  -s\beta B_1 &\cdots  & -s\beta B_q & (\tau+s)I \\
\end{array}\right].
\end{equation}
\end{lemma}

Now, we give a   lemma to   guarantee the positive definiteness of  $G$, defined by
\[
 G=Q+Q\tr -M\tr Q,
\]
 which plays a significant role in showing the whole convergence rate of GS-ADMM.
\begin{lemma} \label{pre-lem-G}
Let $Q, M$ be given by (\ref{Q}) and (\ref{pre-core-M}), respectively. Then,  the matrix
$G$
is symmetric positive definite for any $(\tau,s)\in \mathcal{D}$.
\end{lemma}
\noindent{\bf Proof }
By   simple calculations, the matrix $G$   can be explicitly written as
\[
G=\left[\begin{array}{cc}
H_{\mathbf{x}} & \bf{0} \\
\bf{0}  & \widetilde{G}
\end{array}\right],
\]
 where $H_\mathbf{x}$ is defined in (\ref{Hx}) and
\[
\widetilde{G}=\left[\begin{array}{ccccc|c}
(\sigma_2+1-s)\beta B_1\tr B_1 &&\cdots &&-s\beta B_1\tr B_q &(s-1) B_1\tr \\
   \vdots  &&\ddots &&\vdots&\vdots \\
-s\beta B_q\tr B_1   &&\cdots &&(\sigma_2+1-s)\beta B_q\tr B_q &(s-1)  B_q\tr \\
\hline
     (s-1) B_1 &&\cdots &&(s-1) B_q &\frac{2-\tau-s }{\beta}I
\end{array}\right].
\]
Clearly, the matrix  $G$ is symmetric positive definite if and only if both $H_{\mathbf{x}}$ and $\widetilde{G}$ are symmetric positive definite. Well, $H_{\mathbf{x}}$ is    symmetric and its positivity can be guaranteed by the known conditions that $\sigma_1>p-1$ and the full column rank assumption on the matrices $A_i, i=1,2,\cdots,p$. Hence, we just need to demonstrate the positivity of  the matrix  $\widetilde{G}$.

Noting  that  by the region shown in (\ref{setK})  we have
\begin{equation} \label{regi}
 \tau<1, \quad s<1,\quad \textrm{and}\quad \tau+s<2.
\end{equation}
Besides, it follows
\[
\widetilde{G}=\widetilde{D}\tr \widetilde{G}_0\widetilde{D},
\]
where $\widetilde{D}=\Diag(\beta^{\frac{1}{2}}B_1, \cdots, \beta^{\frac{1}{2}}B_q,\beta^{-\frac{1}{2}}I)$ is a diagonal matrix   and
 \begin{eqnarray*}
\widetilde{G}_0=\left[\begin{array}{ccccccc|c}
(\sigma_2+1-s) I  && -s I  &&\cdots &&-s I  & (s-1)I   \\
-s I  && (\sigma_2+1-s) I  &&\cdots &&-s I  &  (s-1)I   \\
\vdots && \vdots &&\ddots &&\vdots & \vdots  \\
 -s I && -s I  &&\cdots &&(\sigma_2+1-s) I  & (s-1)I   \\
 \hline
 (s-1) I && (s-1) I  &&\cdots &&(s-1) I  & (2-\tau-s)I   \\
\end{array}\right]
 =  P\tr \widetilde{G}_{0,0}P.
 \end{eqnarray*}
 In the above decomposition, we have
 \[
 P=\left[\begin{array}{ccccccc|c}
 I  &&    &&  &&   &     \\
   &&   I  &&  &&  &     \\
  &&   &&\ddots &&  &    \\
  &&    &&  &&  I  &    \\
 \hline
 \frac{1-s}{2-\tau-s} I && \frac{1-s}{2-\tau-s} I  &&\cdots &&\frac{1-s}{2-\tau-s} I  & I  \\
\end{array}\right]
 \]
 and
  \begin{eqnarray*}
\widetilde{G}_{0,0}&=&{\scriptsize\left[\begin{array}{ccccccc|c}
\left(\sigma_2+1-s-\frac{(1-s)^2}{2-\tau-s}\right) I  && \left(-s-\frac{(1-s)^2}{2-\tau-s}\right) I  &&\cdots &&\left(-s-\frac{(1-s)^2}{2-\tau-s}\right)I  & \mathbf{0}   \\
\left(-s-\frac{(1-s)^2}{2-\tau-s}\right) I  && \left(\sigma_2+1-s-\frac{(1-s)^2}{2-\tau-s}\right) I  &&\cdots &&\left(-s-\frac{(1-s)^2}{2-\tau-s}\right) I  &  \mathbf{0}   \\
\vdots && \vdots &&\ddots &&\vdots & \vdots  \\
 \left(-s-\frac{(1-s)^2}{2-\tau-s}\right) I && \left(-s-\frac{(1-s)^2}{2-\tau-s}\right) I  &&\cdots &&\left(\sigma_2+1-s-\frac{(1-s)^2}{2-\tau-s}\right) I  & \mathbf{0}   \\
 \hline
 \mathbf{0} && \mathbf{0}  &&\cdots &&\mathbf{0}  & (2-\tau-s)I   \\
\end{array}\right]}\\
&=& \left[\begin{array}{c|c}
 H_{\mathbf{y},0}+ \left(1-s-\frac{(1-s)^2}{2-\tau-s}\right)EE\tr & \bf{0} \\ \hline
\bf{0}  & (2-\tau-s)I
\end{array}\right],
 \end{eqnarray*}
 where
 \[
E=\left[\begin{array}{c}
I  \\
I \\
\vdots\\
I
\end{array}\right] \quad \mbox{ and } \quad
H_{\mathbf{y},0}=\left[\begin{array}{cccc}
\sigma_2 I  & -I & \cdots &-I \\
-I  & \sigma_2 I  &\cdots &-I \\
 \vdots& \vdots& \ddots& \vdots\\
-I  & -I  &\cdots&  \sigma_2 I
\end{array}\right].
\]
So, the matrix $\widetilde{G}$ is positive definite   if and only if \[H_{\mathbf{y},0}+ \left(1-s-\frac{(1-s)^2}{2-\tau-s}\right)EE\tr\] is positive definite.  Notice that $H_{\mathbf{y},0}$ is positive definite if
 $\sigma_2>q-1$, and  $\left(1-s-\frac{(1-s)^2}{2-\tau-s}\right)EE\tr$ is positive definite if
\[
1-s-\frac{(1-s)^2}{2-\tau-s}=\frac{(1-s)(1-\tau)}{2-\tau-s}>0,
\]
which is clearly guaranteed by the   conditions (\ref{regi}). This completes the proof.
$\ \ \ \diamondsuit$

\begin{theorem} \cite{BaiLiXuZhang2017}\label{Ieq-33}
The sequences  $\{\m{w}^k\}$ and $\{\widetilde{\m{w}}^k\}$ generated by GS-ADMM satisfy
\[
\left\|\m{w}^{k+1}-\m{w}^*\right\|_H^2\leq \left\|\m{w}^{k}-\m{w}^*\right\|_H^2-\left\|\m{w}^{k}-\widetilde{\m{w}}^k\right\|_G^2,\ \forall \m{w}^*\in \mathcal{M}^*,
\]
where $\mathcal{M}^*=\bigcap\limits_{\m{w}\in \mathcal{M}}\left\{\widehat{\m{w}}\in \mathcal{M}|\ h(\m{u})-
h(\widehat{\m{u}})+ \left\langle \m{w}-\widehat{\m{w}}, \mathcal{J}(\m{w})\right\rangle\geq 0\right\}$ and
\begin{equation} \label{HH}
H=QM^{-1}
\end{equation}
is symmetric positive definite for any $(\tau,s)\in \mathcal{D}$.
\end{theorem}

In view of both  Lemma \ref{pre-lem-G} and Theorem  \ref{Ieq-33}, the sequence $\{\widetilde{\m{w}}^k\}$ generated by GS-ADMM is contractive, which implies { a}  global convergence of GS-ADMM. In fact, by estimating the lower bound of $\left\|\m{w}^{k}-\widetilde{\m{w}}^k\right\|_G^2$, { a} global convergence of GS-ADMM {   was} proved   in \cite{BaiLiXuZhang2017}   for the larger region $\C{G}(\supseteq\mathcal{D})$. Next, we will show   { sublinear} nonergodic convergence rate of GS-ADMM for our discussed stepsize region $\mathcal{D}$.

\begin{lemma}  \label{Ieq-33-11}
Let $Q, M, H$ be  given by (\ref{Q}), (\ref{pre-core-M}) and (\ref{HH}), respectively. Then, the sequences  $\{\m{w}^k\}$ and $\{\widetilde{\m{w}}^k\}$ generated by GS-ADMM satisfy
\begin{eqnarray} \label{Contra-wH}
(\m{w}^k-\widetilde{\m{w}}^k)\tr M\tr H M \left\{ (\m{w}^k-\widetilde{\m{w}}^k)-(\m{w}^{k+1}-\widetilde{\m{w}}^{k+1}) \right\}
\geq \frac{1}{2}\left\| (\m{w}^k-\widetilde{\m{w}}^k)-(\m{w}^{k+1}-\widetilde{\m{w}}^{k+1})\right\|_{Q+Q\tr}^2. \nonumber
\end{eqnarray}
\end{lemma}
\noindent{\bf Proof }
 Setting $\m{w}=\widetilde{\m{w}}^{k+1}$ in (\ref{app-vi}), we   obtain
\begin{equation} \label{app-vi-1}
  h(\widetilde{\m{u}}^{k+1})-h(\widetilde{\m{u}}^k)+ \left\langle \widetilde{\m{w}}^{k+1}-\widetilde{\m{w}}^k, \mathcal{J}(\widetilde{\m{w}}^k)+Q(\widetilde{\m{w}}^k-\m{w}^k)\right\rangle \geq 0.
\end{equation}
Meanwhile, the inequality (\ref{app-vi}) with $k:=k+1$ also implies
 \[
  h(\m{u})-h(\widetilde{\m{u}}^{k+1})+ \left\langle \m{w}-\widetilde{\m{w}}^{k+1}, \mathcal{J}(\widetilde{\m{w}}^{k+1})+Q(\widetilde{\m{w}}^{k+1}-\m{w}^{k+1})\right\rangle \geq 0, \quad \forall \m{w}\in \mathcal{M},
 \]
 which, by letting $\m{w}=\widetilde{\m{w}}^{k}$, gives
\begin{equation} \label{app-vi-2}
  h(\widetilde{\m{u}}^{k})-h(\widetilde{\m{u}}^{k+1})+ \left\langle \widetilde{\m{w}}^{k}-\widetilde{\m{w}}^{k+1}, \mathcal{J}(\widetilde{\m{w}}^{k+1})+Q(\widetilde{\m{w}}^{k+1}-\m{w}^{k+1})\right\rangle \geq 0.
\end{equation}
Because of  the skew-symmetric property of $\mathcal{J}(\m{w})$, i.e.,
\[\left\langle \widetilde{\m{w}}^k-\widetilde{\m{w}}^{k+1}, \mathcal{J}(\widetilde{\m{w}}^k)-\mathcal{J}(\widetilde{\m{w}}^{k+1})\right\rangle =0, \]
we have from  (\ref{app-vi-1}) and (\ref{app-vi-2}) that
\begin{equation} \label{app-vi-3}
(\widetilde{\m{w}}^k-\widetilde{\m{w}}^{k+1})\tr Q \left\{ (\m{w}^k-\widetilde{\m{w}}^k)-(\m{w}^{k+1}-\widetilde{\m{w}}^{k+1}) \right\}\geq 0.
\end{equation}
Thus, adding the identity
\begin{eqnarray*}
&&\left\{ (\m{w}^k-\widetilde{\m{w}}^k)-(\m{w}^{k+1}-\widetilde{\m{w}}^{k+1}) \right\}\tr Q \left\{ (\m{w}^k-\widetilde{\m{w}}^k)-(\m{w}^{k+1}-\widetilde{\m{w}}^{k+1}) \right\}\\
&=& \frac{1}{2}\left\| (\m{w}^k-\widetilde{\m{w}}^k)-(\m{w}^{k+1}-\widetilde{\m{w}}^{k+1}) \right\|_{Q+Q\tr}^2
\end{eqnarray*}
to both sides of (\ref{app-vi-3}), we get
\[
 ( \m{w}^k-\m{w}^{k+1} )\tr Q \left\{ (\m{w}^k-\widetilde{\m{w}}^k)-(\m{w}^{k+1}-\widetilde{\m{w}}^{k+1}) \right\} \geq\frac{1}{2}\left\| (\m{w}^k-\widetilde{\m{w}}^k)-(\m{w}^{k+1}-\widetilde{\m{w}}^{k+1}) \right\|_{Q+Q\tr}^2,
\]
which  immediately completes the whole proof by   the relationships   in  (\ref{pre-core}) and (\ref{HH}).
$\ \ \ \diamondsuit$

Next, we establish  the worst-case $\mathcal{O}(1/t)$ nonergodic convergence rate of  GS-ADMM {  in terms of    optimality errors based on the following theorem}.
\begin{theorem} \label{Nonergodic rate}
Let   the sequences  $\{\m{w}^k\}$ and $\{\widetilde{\m{w}}^k\}$ be generated by GS-ADMM. Then, for any integer $t>0$ there exists   a constant $\xi>0$ such that
\[
\left\|M(\m{w}^t-\widetilde{\m{w}}^t)\right\|_H^2\leq \frac{1}{(t+1)\xi}\left\|\m{w}^{0}-\m{w}^*\right\|_H^2,\quad \forall \m{w}^*\in \mathcal{M}^*.
\]
\end{theorem}
\noindent{\bf Proof }
Combining the aforementioned Theorem \ref{Ieq-33} and Lemma \ref{pre-lem-G}, there exists a  constant $\xi>0$ such that
\[
\left\|\m{w}^{k+1}-\m{w}^*\right\|_H^2\leq \left\|\m{w}^{k}-\m{w}^*\right\|_H^2-\xi\left\|M(\m{w}^{k}-\widetilde{\m{w}}^k)\right\|_H^2,\ \forall \m{w}^*\in \mathcal{M}^*,
\]
which suggests
\begin{equation} \label{Rate-121}
\xi\sum\limits_{k=0}^{t}\left\|M(\m{w}^{k}-\widetilde{\m{w}}^k)\right\|_H^2\leq
\left\|\m{w}^{0}-\m{w}^*\right\|_H^2
\end{equation}
for any integer $t>0$.
Meanwhile,  by setting
$
a= M(\m{w}^k-\widetilde{\m{w}}^k) $ and $ b=M(\m{w}^{k+1}-\widetilde{\m{w}}^{k+1})
$
into the following well-known identity
\[
\|a\|_H^2- \|b\|_H^2=2a\tr H(a-b)- \|a-b\|_H^2,
\]
we have
\begin{eqnarray} \label{Ra121}
&&\left\|M(\m{w}^k-\widetilde{\m{w}}^k)\right\|_H^2- \left\|M(\m{w}^{k+1}-\widetilde{\m{w}}^{k+1})\right\|_H^2\nonumber\\
&=&  2 (\m{w}^k-\widetilde{\m{w}}^k)\tr M\tr HM \left\{ (\m{w}^k-\widetilde{\m{w}}^k)-(\m{w}^{k+1}-\widetilde{\m{w}}^{k+1}) \right\} -  \left\| M(\m{w}^k-\widetilde{\m{w}}^k)-M(\m{w}^{k+1}-\widetilde{\m{w}}^{k+1}) \right\|_{H}^2\nonumber\\
&\geq&  \left\| (\m{w}^k-\widetilde{\m{w}}^k)-(\m{w}^{k+1}-\widetilde{\m{w}}^{k+1})\right\|_{Q+Q\tr}^2  -  \left\| M(\m{w}^k-\widetilde{\m{w}}^k)-M(\m{w}^{k+1}-\widetilde{\m{w}}^{k+1}) \right\|_{H}^2\nonumber\\
&=&  \left\| (\m{w}^k-\widetilde{\m{w}}^k)-(\m{w}^{k+1}-\widetilde{\m{w}}^{k+1}) \right\|_{G}^2\nonumber\\
&\geq& 0,
\end{eqnarray}
where
the above first inequality uses Lemma \ref{Ieq-33-11} and  the final equality uses Lemma \ref{pre-lem-G}. Therefore, it holds by (\ref{Ra121}) that
\[
(t+1) \left\|M(\m{w}^{t}-\widetilde{\m{w}}^t)\right\|_H^2 \leq \sum\limits_{k=0}^{t}\left\|M(\m{w}^{k}-\widetilde{\m{w}}^k)\right\|_H^2.
\]
Substituting it into (\ref{Rate-121}), the proof is completed.
$\ \ \ \diamondsuit$

{
\begin{theorem}\label{pointwise-conv-rate}
For any integer $t>0$,  there exists a constant   $\theta>0$  such that
\begin{equation} \label{pointwise-rate}
\left\|\mathbf{d}^t\right\|^2 \le \frac{\theta}{t+1}  \quad \mbox{ and } \quad
\left\| \C{A} \widetilde{\m{x}}^t + \C{B} \widetilde{\m{y}}^t - c\right\|^2 \le  \frac{\theta}{t+1},
\end{equation}
where $\mathbf{d}^t$ is defined by (\ref{right-eq00}) satisfying (\ref{right-141}),
and  $\theta$    depends   on the problem data and the parameters of GS-ADMM.
\end{theorem}
\noindent{\bf Proof }
Let \begin{equation}\label{right-eq00}
\mathbf{d}^t=(\mathbf{d}_{11}^t,\cdots,\mathbf{d}_{1p}^t,\mathbf{d}_{21}^t,
\cdots,\mathbf{d}_{2q}^t),
\end{equation}
componentwisely defined as
\begin{equation}\label{right-eq}
\left\{\begin{array}{ll }
 \mathbf{d}_{1i}^t=-\beta A_i\tr \sum\limits_{l=1,l\neq i}^p A_l
 (\widetilde{x}_l^t-x_l^t )  +\sigma_1\beta A_i\tr\sum\limits_{l=1}^p A_l
 (\widetilde{x}_l^t-x_l^t ),&  i=1,\cdots,p, \\
 \mathbf{d}_{2j}^t=(\sigma_2+1) \beta B_j\tr B_j(\widetilde{y}_j^{t}-y_j^t)-\tau B_j\tr (\widetilde{\lambda}^t-\lambda^t),& j=1,\cdots,q.
\end{array}\right.
\end{equation}
Then, according to the proof of \cite[Lemma 2]{BaiLiXuZhang2017},  that is, the first-order optimality conditions of  the  subproblems of GS-ADMM, we have
\[
\left\{\begin{array}{ll }
  f_i(x_i)-f_i(\widetilde{x}_i^t)
+(x_i-\widetilde{x}_i^t)\tr (-A_i\tr  \widetilde{\lambda}^t + \mathbf{d}_{1i}^t)\geq 0,&  i=1,\cdots,p, \\
  g_j(y_j)-g_j(\widetilde{y}_j^{t}) +(y_j-\widetilde{y}_j^{t})\tr (-B_j\tr  \widetilde{\lambda}^t + \mathbf{d}_{2j}^t)\geq 0,& j=1,\cdots,q,
\end{array}\right.
\]
which   implies
\begin{equation}\label{right-141}
 A_i\tr  \widetilde{\lambda}^t - \mathbf{d}_{1i}^t\in \partial f_i(\widetilde{x}_i^t)+\mathcal{N}_{\mathcal{X}_i}(\widetilde{x}_i^t)\quad \textrm{and} \quad
B_j\tr  \widetilde{\lambda}^t - \mathbf{d}_{2j}^t\in \partial g_j(\widetilde{y}_j^{t}) + \mathcal{N}_{\mathcal{Y}_j}(\widetilde{y}_j^{t}).
\end{equation}
Here the notation $\mathcal{N}_{\mathcal{X}}(x)$ denotes the normal cone of $\mathcal{X}$ at $x$.
By (\ref{right-eq}) and Theorem \ref{Nonergodic rate}, it can be deduced  that
\[
\left\|\mathbf{d}^t\right\|^2= \sum\limits_{i=1}^{p}\left\|\mathbf{d}_{1i}^t\right\|^2 +\sum\limits_{j=1}^{q}\left\|\mathbf{d}_{2j}^t\right\|^2\leq \theta \left\|\widetilde{\m{w}}^t-\m{w}^{t}\right\|^2\leq \frac{\theta}{t+1},
\]
where and in the following proof, $\theta$    depends only on the problem data and the parameters of GS-ADMM.

  We next prove the inequality in the right-hand side of (\ref{pointwise-rate}). Since the equality (\ref{tilde-lambda}) can be rewritten as
\[
\left(\mathcal{A}\widetilde{\m{x}}^{k}+\mathcal{B}\widetilde{\m{y}}^{k}-c\right)-
\mathcal{B}\left(\widetilde{\m{y}}^{k}-\m{y}^k\right)
+\frac{1}{\beta} (\widetilde{\lambda}^k-\lambda^k)=0,
\]
 we have
\[
\left\| \C{A} \widetilde{\m{x}}^t + \C{B} \widetilde{\m{y}}^t - c\right\|^2
=\left\|\mathcal{B}\left(\widetilde{\m{y}}^{t}-\m{y}^t\right)
-\frac{1}{\beta} (\widetilde{\lambda}^t-\lambda^t)\right\|^2 \le \frac{\theta}{t+1}.\ \ \ \diamondsuit
\]
}

Clearly, a nonergodic convergence rate in general is stronger  { than} the ergodic convergence rate for GS-ADMM.  Let $c_0=\inf\limits_{\m{w}^*\in \mathcal{M}^*}\left\|\m{w}^{0}-\m{w}^*\right\|_H^2.$ Then, for any tolerance $\epsilon>0$, Theorem  \ref{Nonergodic rate} tells us that   it needs at most $\lfloor \frac{c_0}{\xi\epsilon}\rfloor$ iterations to ensure $\left\|M(\m{w}^t-\widetilde{\m{w}}^t)\right\|_H^2\leq \epsilon.$ { If $\C{A} \widetilde{\m{x}}^t + \C{B} \widetilde{\m{y}}^t = c$ and $\mathbf{d}^t=0$, then we will have $\widetilde{\m{w}}^t \in \mathcal{M}^*$. Hence, we could use $\widetilde{\m{w}}^t$ (or equivalently the iterate $\m{w}^t$ since $\lim\limits_{t\rightarrow\infty}\widetilde{\m{w}}^t-\m{w}^t=0$ by the proof of \cite[Theorem 6]{BaiLiXuZhang2017}) as an approximate solution
of the problem when the
right-hand sides of the inequalities in (\ref{pointwise-rate}) are sufficiently small. }

\subsection{Linear convergence rate}
Throughout this subsection, all  subdifferentials of the functions $f_i, g_j$ in (\ref{Prob}) are assumed to be piecewise liner multi-functions. Under this hypothesis we will prove { a} global linear convergence rate of GS-ADMM by the aid of an error function
\[
\textrm{dist}^2_H(\m{w}^k, \mathcal{M}^*):=\min\left\{\left\|\m{w}^k-\m{w}^*\right\|^2_H\ | \ \m{w}^*\in\mathcal{M}^*\right\}.
\]
If $H=I,$ we simply denote $\textrm{dist}^2_{I}(\m{w}^k, \mathcal{M}^*)$ by $\textrm{dist}^2(\m{w}^k, \mathcal{M}^*).$

Since each   $\mathcal{X}_i$ in the problem (\ref{Prob}) is   a polyhedron,  so $\mathcal{X}_i$ is convex   and any projection operator $\mathcal{P}_{\mathcal{X}_i}(x_i):=\arg\min\limits_{c\in\mathcal{X}_i}\|c-x_i\| $
 is piecewise linear from  { \cite[Proposition 4.1.4]{FPang2003}}. Here
  $\mathcal{P}_{\mathcal{X}_i}$ is nonexpansive, that is, the following inequality holds:
\[
\left\|\mathcal{P}_{\mathcal{X}_i}(c)-\mathcal{P}_{\mathcal{X}_i}(x_i)\right\|\leq \|c-x_i\|, \quad \forall c,x_i\in\mathbb{R}^{m_i}.
\]

Let $\partial f(x)$
be  the sub-differential
of a convex function $f(x):\mathbb{R}^n\rightarrow \mathbb{R}$, defined as
\[
\partial f(x):=\left\{ d\in \mathbb{R}^{n}|\  f(v)\geq f(x)+ d\tr (v-x),\ \forall  v\in \mathbb{R}^{n}\right\}.
\]
Then, for any saddle-point $\m{w}^*=(x_i^*,\cdots,x_p^*,y_1^*,\cdots,y_q^*,\lambda^*)$ of   (\ref{Prob}), { there} exist $\eta_i\in\partial f_i(x_i) (i=1,\cdots,p)$ and $ \nu_j\in\partial g_j(y_j) (j=1,\cdots,q)$ such that
\[
\left(\begin{array}{c}
x_1-x_1^*\\ \vdots\\ x_p-x_p^*\\y_1-y_1^*\\ \vdots\\ y_q-y_q^*\\ \lambda-\lambda^*
\end{array}\right)\tr
\left(\begin{array}{c}
\eta_1-A_1\tr \lambda^*\\ \vdots\\ \eta_p-A_p\tr \lambda^*\\\nu_1-B_1\tr \lambda^*\\ \vdots\\ \nu_q-B_q\tr \lambda^*\\ \mathcal{A} \mathbf{x}^*+\mathcal{B}\mathbf{y}^*-c
\end{array}\right) \geq 0,\quad\forall \m{w}\in \mathcal{M},
\]
which can be characterized by solving the equation $ \left\|e_{\mathcal{M}}(\m{w},\gamma)\right\|=0$
with
\[
e_{\mathcal{M}}(\m{w},\gamma)=\left(\begin{array}{ccc}
e_{\mathcal{X}_1}(\m{w},\gamma):&=& x_1- \mathcal{P}_{\mathcal{X}_1}\left[x_1 -\gamma (\partial f_1(x_1)-A_1\tr \lambda)\right]\\
 \vdots&&\\
e_{\mathcal{X}_p}(\m{w},\gamma):&=& x_p- \mathcal{P}_{\mathcal{X}_p}\left[x_p -\gamma (\partial f_p(x_p)-A_p\tr \lambda)\right]\\
e_{\mathcal{Y}_1}(\m{w},\gamma):&=& y_1- \mathcal{P}_{\mathcal{Y}_1}\left[y_1 -\gamma (\partial g_1(y_1)-B_1\tr \lambda)\right]\\
 \vdots&&\\
e_{\mathcal{Y}_q}(\m{w},\gamma):&=& y_q- \mathcal{P}_{\mathcal{Y}_q}\left[y_q -\gamma (\partial g_q(y_q)-B_q\tr \lambda)\right]\\
 e_{\Lambda}(\m{w},\gamma): &=&\gamma \left(\mathcal{A} \mathbf{x}+\mathcal{B}\mathbf{y}-c\right)
\end{array}\right),\ \forall \gamma>0.
\]
Under the assumption that
$\partial f_i$ and $ \partial g_j$ are piecewise linear multi-functions, $e_{\mathcal{M}}(\m{w},\gamma)$ is  also   piecewise linear. Besides,
$\m{w}^* \in \mathcal{M}^*$ if and only if $ \textbf{0} \in e_{\mathcal{M}}(\m{w},1)$.
The following lemma,   coming from Robinsons's continuity property \cite{Robinson1981}
for polyhedral multi-functions,
shows that
$\textrm{dist}(\mathbf{0},e_{\mathcal{M}}(\m{w},1))$ could  provide a global
error bound on the distance of $\m{w}$ to the   solution set $\mathcal{M}^*$.
\begin{lemma} \label{bjc}
Under the assumption that
$\partial f_i$ and $ \partial g_j$ are piecewise linear multi-functions, there exists a constant $\zeta>0$ such that
\[
\textrm{dist}(\m{w},\mathcal{M}^*)\leq \zeta \textrm{dist}(\mathbf{0},e_{\mathcal{M}}(\m{w},1)), \ \forall \m{w}\in\mathcal{M}.
\]
\end{lemma}

For   convenience of analysis, let
\begin{equation}\label{S2-000}
\widetilde{\mu}_i=\lambda_{\max}(A_i\tr A_i)~\ \textrm{and}~\  \widetilde{\nu}_j=\lambda_{\max}(B_j\tr B_j), \quad \forall i=1,2,\cdots,p,     j=1,2,\cdots,q.
\end{equation}
Define
\begin{equation}\label{S2-001}
\delta:=\max\left\{\ \overline{\theta}_i, \overline{\vartheta}_j, \overline{\eta}\ | i=1,2,\cdots,p,    j=1,2,\cdots,q\right\}
\end{equation}
with
\[
  \left \{\begin{array}{lll}
\overline{\theta}_i&=& 4p(1-\sigma_1)^2\beta^2 \sum\limits_{l=1}^{p} \widetilde{\mu}_l+ 4\beta^2 \widetilde{\mu}_i, \quad \forall i=1,2,\cdots,p, \\
\overline{\vartheta}_j &=&4q(s\beta)^2\sum\limits_{l=1}^{p} \widetilde{\mu}_l
 +3q(s\beta)^2\sum\limits_{l=1}^{q}\widetilde{\nu}_l
 + 3(\sigma_2+1)^2\beta^2 \widetilde{\nu}_j+2q, \quad \forall j=1,2,\cdots,q,  \\
 \overline{\eta}&=& 4(\tau+s-1)^2\sum\limits_{l=1}^{p} \widetilde{\mu}_l
   +3(s-1)^2\sum\limits_{l=1}^{q}\widetilde{\nu}_l+\frac{2}{\beta^2}.
\end{array}\right.
\]
Note that all the above $(p+q+1)$ notations are positive since   the matrices $A_i, B_j$ have full column rank. Hence,   $\delta$ is a positive number.
\begin{theorem}\label{rate-linear}
Let  $\delta$ be defined in (\ref{S2-001}) with $\widetilde{\mu}_i, \widetilde{\nu}_j$ being defined in (\ref{S2-000}). Then, the sequences $\{\m{w}^k\}$ and $\{\widetilde{\m{w}}^k\}$  generated by GS-ADMM satisfy
\[
\textrm{dist}^2(\mathbf{0},e_{\mathcal{M}}(\m{w}^{k+1},1))\leq \frac{  \max\limits_{i,j}\left\{\delta\widetilde{\mu}_i,\delta\widetilde{\nu}_j,\delta\right\}}
{\lambda_{\min}(G)}\left\|\m{w}^k-\widetilde{\m{w}}^k\right\|_{G}^2,\quad \forall i=1,\cdots,p,    j=1,\cdots,q.
\]
\end{theorem}
\noindent{\bf Proof }
Firstly, by the equation (20) mentioned in \cite{BaiLiXuZhang2017}, that is,
\[
f_i(x_i)-f_i(\widetilde{x}_i^k)
+ \left\langle A_i (x_i-\widetilde{x}_i^k),
 - \widetilde{\lambda}^k  -\beta \sum\limits_{l=1,l\neq i}^p A_l
 (\widetilde{x}_l^k-x_l^k )  +\sigma_1\beta \sum\limits_{l=1}^p A_l
 (\widetilde{x}_l^k-x_l^k )  \right\rangle\geq 0,\ \forall x_i\in\mathcal{X}_i,
\]
there exists $\eta_i\in\partial f_i(x_i), i=1,\cdots,p,$ such that
\[
x_i^{k+1}=\mathcal{P}_{\mathcal{X}_i}\left[x_i^{k+1}-\left(\eta_i-A_i\tr \widetilde{\lambda}^k-\beta A_i\tr\sum\limits_{l=1,l\neq i}^p A_l
 (\widetilde{x}_l^k-x_l^k )  +\sigma_1\beta A_i\tr\sum\limits_{l=1}^p A_l
 (\widetilde{x}_l^k-x_l^k )\right) \right].
\]
Therefore, we have from the definition of $\textrm{dist}(\mathbf{0},\cdot)$ and the nonexpansive property of the projection operator that
\begin{eqnarray} \label{S2-003}
&&\textrm{dist}(\mathbf{0},e_{\mathcal{X}_i}(\m{w}^{k+1},1)) \nonumber\\
& =& \textrm{dist}\left(x_i^{k+1}, \mathcal{P}_{\mathcal{X}_i}\left[x_i^{k+1}-\left\{\partial f_i(x_i^{k+1})-A_i\tr \lambda^{k+1}\right\}\right]\right) \nonumber\\
 &\leq& \left\|
 \begin{array}{l}
 \mathcal{P}_{\mathcal{X}_i}\left[x_i^{k+1}-\left(\eta_i-A_i\tr \widetilde{\lambda}^k-\beta A_i\tr\sum\limits_{l=1,l\neq i}^p A_l
 (\widetilde{x}_l^k-x_l^k )  +\sigma_1\beta A_i\tr\sum\limits_{l=1}^p A_l
 (\widetilde{x}_l^k-x_l^k )\right) \right]\\
\quad -\mathcal{P}_{\mathcal{X}_i}\left[x_i^{k+1}-\left\{\eta_i-A_i\tr \lambda^{k+1}\right\}\right]
 \end{array}
 \right\|\nonumber\\
&\leq&  \left\|  A_i\tr (\widetilde{\lambda}^k- \lambda^{k+1})+\beta A_i\tr\sum\limits_{l=1,l\neq i}^p A_l
 (\widetilde{x}_l^k-x_l^k )  -\sigma_1\beta A_i\tr\sum\limits_{l=1}^p A_l
 (\widetilde{x}_l^k-x_l^k )\right\|\nonumber\\
&= & \left\|  \begin{array}{l}
A_i\tr \left[(\tau+s-1)( \lambda^{k}-\widetilde{\lambda}^k)-s\beta \sum\limits_{j=1}^{q}B_j(y_j^k-\widetilde{y}_j^k)\right]\\
\quad  + A_i\tr \left[ (1-\sigma_1)\beta  \sum\limits_{l=1}^p A_l
 (\widetilde{x}_l^k-x_l^k )-\beta   A_i(\widetilde{x}_i^k-x_i^k )\right]
 \end{array}\right\|,
\end{eqnarray}
where the second equality uses the fact
\begin{eqnarray} \label{S2-004}
\begin{array}{lll}
\lambda^{k+1}  &=& \lambda^{k+\frac{1}{2}}-s\beta(\mathcal{A}\m{x}^{k+1}+\mathcal{B}\m{y}^{k+1}-c)    \\
 &=&  \lambda^{k+\frac{1}{2}}-s\beta(\mathcal{A}\m{x}^{k+1}+\mathcal{B}\m{y}^{k}-c)+s\beta \mathcal{B}(\m{y}^{k}-\m{y}^{k+1})  \\
 &=& \lambda^k-\tau(\lambda^k- \widetilde{\lambda}^k)-s(\lambda^k- \widetilde{\lambda}^k)+s\beta \mathcal{B}(\m{y}^{k}-\widetilde{\m{y}}^{k}) \\
&=& \lambda^k-\left[-s\beta \sum\limits_{j=1}^{q}B_j(y_j^k-\widetilde{y}_j^k)+(\tau+s)(\lambda^k- \widetilde{\lambda}^k) \right].
\end{array}
\end{eqnarray}
Similarly, there exists $ \nu_j\in\partial g_j(y_j),j=1,\cdots,q,$ such that
\[
y_j^{k+1}=\mathcal{P}_{\mathcal{Y}_j}\left[y_j^{k+1}-\left(\nu_j-B_j\tr\widetilde{\lambda}^k+(\sigma_2+1)\beta B_j\tr B_j
 (\widetilde{y}_j^k-y_j^k )-\tau B_j\tr(\widetilde{\lambda}^k-\lambda^k)\right) \right].
\]
Hence, we have
\begin{eqnarray} \label{Sec33-101}
&&\textrm{dist}(\mathbf{0},e_{\mathcal{Y}_j}(\m{w}^{k+1},1)) \nonumber\\
 &\leq& \left\|
 \begin{array}{l}
 \mathcal{P}_{\mathcal{Y}_j}\left[y_j^{k+1}-\left(\nu_j-B_j\tr\widetilde{\lambda}^k+(\sigma_2+1)\beta B_j\tr B_j
 (\widetilde{y}_j^k-y_j^k )-\tau B_j\tr(\widetilde{\lambda}^k-\lambda^k)\right) \right]\\
\quad -\mathcal{P}_{\mathcal{Y}_j}\left[y_j^{k+1}-\left\{\nu_j-B_j\tr \lambda^{k+1}\right\}\right]
 \end{array}
 \right\|\nonumber\\
&\leq&  \left\|  B_j\tr (\widetilde{\lambda}^k- \lambda^{k+1})-(\sigma_2+1)\beta B_j\tr B_j
 (\widetilde{y}_j^k-y_j^k )+\tau B_j\tr(\widetilde{\lambda}^k-\lambda^k) \right\|\nonumber\\
&= &
\left\| \begin{array}{l}
B_j\tr \left[(\tau+s-1)( \lambda^{k}-\widetilde{\lambda}^k)-s\beta \sum\limits_{j=1}^{q}B_j(y_j^k-\widetilde{y}_j^k)\right]\\
\quad -(\sigma_2+1)\beta B_j\tr B_j
 (\widetilde{y}_j^k-y_j^k )+\tau B_j\tr(\widetilde{\lambda}^k-\lambda^k)
\end{array} \right\|\nonumber\\
 &= &
\left\|
B_j\tr \left[(s-1)( \lambda^{k}-\widetilde{\lambda}^k)-s\beta \sum\limits_{j=1}^{q}B_j(y_j^k-\widetilde{y}_j^k)-(\sigma_2+1)\beta  B_j
 (\widetilde{y}_j^k-y_j^k )\right]
 \right\|.
\end{eqnarray}

Secondly,  we can get by the update of $\lambda^{k+1}$ and $\lambda^{k+\frac{1}{2}}$ in GS-ADMM as well as (\ref{tilde-lambda}) and  (\ref{S2-004}) that
\[
\begin{array}{lll}
&&\mathcal{A}\m{x}^{k+1}+\mathcal{B}\m{y}^{k+1}-c   = \frac{1}{s\beta}\left(\lambda^{k+\frac{1}{2}}-\lambda^{k+1}\right)    \\
 &=&  \frac{1}{s\beta}\left[\lambda^k-\lambda^{k+1}-\tau\beta\left(\mathcal{A}\m{x}^{k+1}+\mathcal{B}\m{y}^{k}-c\right)\right]   \\
 &=& \frac{1}{s\beta}\left[-s\beta \sum\limits_{j=1}^{q}B_j(y_j^k-\widetilde{y}_j^k)+(\tau+s)(\lambda^k- \widetilde{\lambda}^k) -\tau(\lambda^k- \widetilde{\lambda}^k)\right]\\
&=& \frac{1}{\beta}(\lambda^k- \widetilde{\lambda}^k)- \sum\limits_{j=1}^{q}B_j(y_j^k-\widetilde{y}_j^k),
\end{array}
\]
which further shows
\begin{eqnarray} \label{S2-006}
\textrm{dist}(\mathbf{0},e_{\Lambda}(\m{w}^{k+1},1)) = \left\|\mathcal{A}\m{x}^{k+1}+\mathcal{B}\m{y}^{k+1}-c\right\|=  \left\|\frac{1}{\beta}(\lambda^k- \widetilde{\lambda}^k)- \sum\limits_{j=1}^{q}B_j(y_j^k-\widetilde{y}_j^k)  \right\|.
\end{eqnarray}

{ Denote}  by
\begin{eqnarray} \label{S2-106}
   \left\{\begin{array}{llll}
   \widetilde{z}^k   &: = & \lambda^k- \widetilde{\lambda}^k,&\\
 \widetilde{a}_i^k&: = &x_i^k-\widetilde{x}_i^k, &   i=1,2,\cdots,p,\\
\widetilde{b}_j^k& := & y_j^k-\widetilde{y}_j^k,& j=1,2,\cdots,q.
\end{array}\right.
\end{eqnarray}
Then, by combining (\ref{S2-003}), (\ref{Sec33-101})-(\ref{S2-106}) together with the following identity
\[
\left\|d_1+d_2+\cdots+d_n\right\|^2\leq n\left(\|d_1\|^2+\|d_2\|^2+\cdots+\|d_n\|^2\right)
\]
for any $d_i\in \mathbb{R}^n, i=1,2\cdots,,n$, it can be achieved by the fact $\lambda_{\max(A\tr A)}= \lambda_{\max(AA\tr)}$ that
 \begin{eqnarray*}
&&\textrm{dist}^2(\mathbf{0},e_{\mathcal{M}}(\m{w}^{k+1},1)) \\
& =& \sum\limits_{i=1}^{p}\textrm{dist}^2(\mathbf{0},e_{\mathcal{X}_i}(\m{w}^{k+1},1))+
     \sum\limits_{j=1}^{q}\textrm{dist}^2(\mathbf{0},e_{\mathcal{Y}_j}(\m{w}^{k+1},1))+
     \textrm{dist}^2(\mathbf{0},e_{\Lambda}(\m{w}^{k+1},1)) \\
&\leq & \sum\limits_{i=1}^{p} \widetilde{\mu}_i\left\| (\tau+s-1)\widetilde{z}^k
    -s\beta \sum\limits_{l=1}^{q}B_l\widetilde{b}_l^k
    -\beta  A_i\widetilde{a}_i^k +
    (1-\sigma_1)\beta \sum\limits_{l=1}^p A_l \widetilde{a}_l^k\right\|^2  \\
 &&+\sum\limits_{j=1}^{q}\widetilde{\nu}_j \left\|(s-1)\widetilde{z}^k
 -s\beta \sum\limits_{l=1}^{q}B_l\widetilde{b}_l^k-
 (\sigma_2+1)\beta  B_j\widetilde{b}_j^k\right\|^2
 +\left\|\frac{1}{\beta}\widetilde{z}^k- \sum\limits_{j=1}^{q}B_j\widetilde{b}_j^k  \right\|^2 \\
&\leq & 4\sum\limits_{i=1}^{p} \widetilde{\mu}_i\left\{ (\tau+s-1)^2\left\|\widetilde{z}^k\right\|^2
    +(s\beta)^2 \left\|\sum\limits_{l=1}^{q}B_l\widetilde{b}_l^k\right\|^2
    +\beta^2\left\|A_i\widetilde{a}_i^k\right\|^2
    +(1-\sigma_1)^2\beta^2 \left\|\sum\limits_{l=1}^p A_l \widetilde{a}_l^k\right\|^2 \right\}  \\
&&+ 3\sum\limits_{j=1}^{q}\widetilde{\nu}_j \left\{(s-1)^2\left\|\widetilde{z}^k\right\|^2
 +(s\beta)^2 \left\|\sum\limits_{l=1}^{q}B_l\widetilde{b}_l^k\right\|^2
 +
 (\sigma_2+1)^2\beta^2  \left\|B_j\widetilde{b}_j^k\right\|^2\right\} \\
 &&+2 \left\{ \frac{1}{\beta^2}\left\| \widetilde{z}^k\right\|^2+ \left\|\sum\limits_{j=1}^{q}B_j\widetilde{b}_j^k \right\|^2\right\}
\end{eqnarray*}
\begin{eqnarray*}
&\leq & 4\sum\limits_{i=1}^{p} \widetilde{\mu}_i
     \left\{ (\tau+s-1)^2\left\|\widetilde{z}^k\right\|^2
    +q(s\beta)^2 \sum\limits_{l=1}^{q}\left\|B_l\widetilde{b}_l^k\right\|^2
    +\beta^2\left\| A_i\widetilde{a}_i^k\right\|^2
    +p(1-\sigma_1)^2\beta^2 \sum\limits_{l=1}^p\left\| A_l \widetilde{a}_l^k\right\|^2 \right\}  \\
&&+ 3\sum\limits_{j=1}^{q}\widetilde{\nu}_j\left\{(s-1)^2\left\|\widetilde{z}^k\right\|^2
 +q(s\beta)^2 \sum\limits_{l=1}^{q}\left\|B_l\widetilde{b}_l^k\right\|^2
 +
 (\sigma_2+1)^2\beta^2  \left\|B_j\widetilde{b}_j^k\right\|^2\right\} \\
 &&+2 \left\{ \frac{1}{\beta^2}\left\| \widetilde{z}^k\right\|^2+ q\sum\limits_{j=1}^{q}\left\|B_j\widetilde{b}_j^k \right\|^2\right\} \\
&=&  \sum\limits_{j=1}^{q}\left\{
 4q(s\beta)^2\sum\limits_{l=1}^{p} \widetilde{\mu}_l
 +3q(s\beta)^2\sum\limits_{l=1}^{q}\widetilde{\nu}_l
 + 3(\sigma_2+1)^2\beta^2 \widetilde{\nu}_j
 +2q
 \right\}\left\|B_j\widetilde{b}_j^k\right\|^2  \\
&&+
   \left\{ 4(\tau+s-1)^2\sum\limits_{l=1}^{p} \widetilde{\mu}_l
   +3(s-1)^2\sum\limits_{l=1}^{q}\widetilde{\nu}_l+\frac{2}{\beta^2}\right\} \left\|\widetilde{z}^k\right\|^2  \\
&&+  \sum\limits_{i=1}^{p}\left\{ 4p(1-\sigma_1)^2\beta^2 \sum\limits_{l=1}^{p} \widetilde{\mu}_l+ 4\beta^2 \widetilde{\mu}_i    \right\}
  \left\| A_i \widetilde{a}_i^k\right\|^2  \\
&=&\sum\limits_{i=1}^{p} \overline{\theta}_i\left\| A_i \widetilde{a}_i^k\right\|^2+
  \sum\limits_{j=1}^{q} \overline{\vartheta}_j\left\|B_j\widetilde{b}_j^k\right\|^2+
  \overline{\eta}\left\|\widetilde{z}^k\right\|^2  \\
&\leq &\delta \max\limits_{i,j}\left\{\widetilde{\mu}_i,\widetilde{\nu}_j,1\right\}
\left\|\m{w}^k-\widetilde{\m{w}}^k\right\|^2
\leq  \frac{ \max\limits_{i,j}\left\{\delta\widetilde{\mu}_i,\delta\widetilde{\nu}_j,\delta\right\}}
{\lambda_{\min}(G)}\left\|\m{w}^k-\widetilde{\m{w}}^k\right\|_{G}^2. \ \ \ \diamondsuit
\end{eqnarray*}

Based on the above preparations, we   show   { a} global linear convergence rate of GS-ADMM.
\begin{theorem}\label{rate-linear-11}
Let $\delta$ be defined in (\ref{S2-001}) with $\widetilde{\mu}_i, \widetilde{\nu}_j$ being defined in (\ref{S2-000}). Then,
there exists a constant $\zeta>0$ such that  the   sequence $\{\m{w}^{k}\}$   generated by GS-ADMM satisfies
\[
\textrm{dist}^2_{H}(\m{w}^{k+1}, \mathcal{M}^*)\leq \frac{1}{1+\epsilon }\textrm{dist}^2_{H}(\m{w}^k, \mathcal{M}^*),
\]
where
\[
\epsilon= \frac{\lambda_{\min}(G)}{\zeta^2 \max\limits_{i,j}\left\{\delta\widetilde{\mu}_i,\delta\widetilde{\nu}_j,\delta\right\}\lambda_{\max}(H)}>0, \quad \forall i=1,\cdots,p,    j=1,\cdots,q.
\]
\end{theorem}
\noindent{\bf Proof }
Because $ \mathcal{M}^*$ is a closed convex  set, there exists a $\m{w}^*_k\in \mathcal{M}^*$ satisfying  \[\textrm{dist}_{H}(\m{w}^k, \mathcal{M}^*)=\left\|\m{w}^k-\m{w}^*_k\right\|_H.\] Then,
by Lemma \ref{bjc} and Theorem \ref{rate-linear} there exists a constant $\zeta>0$ such that
\begin{eqnarray} \label{Sec33-108}
\left\|\m{w}^k-\widetilde{\m{w}}^k\right\|_{G}
&\geq&  \sqrt{  \frac{\lambda_{\min}(G)}{ \max\limits_{i,j}\left\{\delta\widetilde{\mu}_i,\delta\widetilde{\nu}_j,\delta\right\}} } \textrm{dist}(\mathbf{0},e_{\mathcal{M}}(\m{w}^{k+1},1))\nonumber\\
&\geq&  \sqrt{  \frac{\lambda_{\min}(G)}{\zeta^2  \max\limits_{i,j}\left\{\delta\widetilde{\mu}_i,\delta\widetilde{\nu}_j,\delta\right\}} } \textrm{dist}(\m{w}^{k+1} , \mathcal{M}^*)\nonumber\\
&\geq&  \sqrt{  \frac{\lambda_{\min}(G)}{\zeta^2 \max\limits_{i,j}\left\{\delta\widetilde{\mu}_i,\delta\widetilde{\nu}_j,\delta\right\}\lambda_{\max}(H)} } \textrm{dist}_H(\m{w}^{k+1} , \mathcal{M}^*), \nonumber
\end{eqnarray}
where $G$ and $H$ are respectively defined in Lemma \ref{pre-lem-G} and (\ref{HH}). So, we will have from the above inequality that
\begin{eqnarray*} \label{Sec33-1090}
&&(1+\epsilon )\textrm{dist}^2_{H}(\m{w}^{k+1}, \mathcal{M}^*) \nonumber\\
&\leq&  \left\|\m{w}^{k+1}-\m{w}^*_k\right\|_H^2 + \epsilon \textrm{dist}^2_{H}(\m{w}^{k+1}, \mathcal{M}^*) \\
&\leq& \|\m{w}^{k+1}-\m{w}^*_k\|_H^2 + \epsilon   \frac{\zeta^2 \max\limits_{i,j}\left\{\delta\widetilde{\mu}_i,\delta\widetilde{\nu}_j,\delta\right\}\lambda_{\max}(H)}{\lambda_{\min}(G)}  \left\|\m{w}^k-\widetilde{\m{w}}^k\right\|_{G}^2 \\
&=& \left\|\m{w}^{k+1}-\m{w}^*_k\right\|_H^2 + \|\m{w}^k-\widetilde{\m{w}}^k\|_{G}^2 \\
&\leq&\left\|\m{w}^{k}-\m{w}^*_k\right\|_H^2 = \textrm{dist}^2_{H}(\m{w}^k, \mathcal{M}^*).
\end{eqnarray*}
This completes the whole proof. $  \ \ \ \diamondsuit$

 Next, we show that $\{\m{w}^k\}$ generated by GS-ADMM converges to a point
$\m{w}^\infty \in \mathcal{M}^*$ R-linearly.
\begin{corollary} \label{13131}
  Let $\epsilon >0$ be  defined in Theorem \ref{rate-linear-11} and the  sequence   $\{\m{w}^k\}$ be generated by
GS-ADMM. Then, there exists a point $\m{w}^\infty \in \mathcal{M}^*$
such that
\begin{equation}\label{v-lin-conv}
\left\|\m{w}^k-\m{w}^\infty \right\|_{H}\leq C r^k,
\end{equation}
where
\[
C=\frac{2\textrm{dist}_{H}(\m{w}^0, \mathcal{M}^*)}{1-r}>0\quad \mbox{and}\quad r=\frac{1}{\sqrt{1+\epsilon}}\in(0,1).
\]
\end{corollary}
\noindent{\bf Proof }
Select $\m{w}_k^*\in \mathcal{M}^*$ such that $\textrm{dist}_{H}(\m{w}^k, \mathcal{M}^*)=\left\|\m{w}^k-\m{w}_k^*\right\|_{H}$ and let
\begin{equation}\label{abbc}
\m{w}^{k+1}=\m{w}^{k}+d^k.
\end{equation}
Then, it follows from Theorem \ref{Ieq-33} that
$
\left\|\m{w}^{k+1}-\m{w}_k^*\right\|_{H}\leq \left\|\m{w}^k-\m{w}_k^*\right\|_{H}
$
  implying
\begin{eqnarray}\label{d-linear}
\left\|d^k\right\|_{\widetilde{Q}}&=&\left\|\m{w}^{k+1}-\m{w}^k\right\|_{H} \nonumber\\
&\leq & \left\|\m{w}^{k+1}-\m{w}_k^*\right\|_{H}+ \left\|\m{w}^k-\m{w}_k^*\right\|_{H} \nonumber\\
&\leq & 2\left\|\m{w}^k-\m{w}_k^*\right\|_{H}
= 2\textrm{dist}_{H}(\m{w}^k, \mathcal{M}^*)\nonumber\\
&\leq &2 r^k \textrm{dist}_{H}(\m{w}^0, \mathcal{M}^*),
\end{eqnarray}
where the last inequality comes  from Theorem \ref{rate-linear-11}.
According to  { \cite[Theorem 6]{BaiLiXuZhang2017}}, the sequence  $\{\m{w}^k\} $
generated by GS-ADMM
converges to a $ w^\infty  \in \mathcal{M}^*$.
Hence, we obtain by (\ref{abbc}) that $ \m{w}^\infty =\m{w}^k+\sum_{j=k}^{\infty}d^j $,
which together with (\ref{d-linear}) show
\begin{eqnarray*}
\left\|\m{w}^k-\m{w}^\infty\right\|_{H}&\leq&\sum\limits_{j=k}^{\infty}\|d^j\|_{H}
\leq  2\textrm{dist}_{H}(\m{w}^0, \mathcal{M}^*)\sum\limits_{j=k}^{\infty}r^k  \\
&= &  2\textrm{dist}_{H}(\m{w}^0, \mathcal{M}^*) r^k \sum\limits_{j=0}^{\infty}r^k\\
 &\leq&  r^k \left[2\textrm{dist}_{H}(\m{w}^0, \mathcal{M}^*)\frac{1}{1-r}\right].
\end{eqnarray*}
Hence, the assertion (\ref{v-lin-conv}) holds, namely, $\m{w}^k$ converges $\m{w}^\infty$ R-linearly.  $\ \ \ \diamondsuit$

\section{Conclusion remark}

In this note, we further study   iteration-complexity of   GS-ADMM for solving the prototype multi-block separable convex optimization model. We establish its   sublinear nonergodic convergence rate and also a R-linear convergence rate under   assumptions that the sub-differential of each component function in the objective
function is piecewise linear and all the constraint sets are polyhedra.
By the fourth part discussed in \cite{BaiLiXuZhang2017} and the  analysis in this work, the GS-ADMM with either $\sigma_1=0$ or $\sigma_2=0$ has a similar convergence rate as described in Theorem \ref{pointwise-conv-rate}, Theorem \ref{rate-linear-11} and Corollary \ref{13131}.
Viewed from the proof of Theorem \ref{rate-linear-11},  the  linear convergence analysis depends mainly on  Theorem \ref{rate-linear} and the positivity of the matrix $G$. Hence,   if the sequence generated by an algorithm has the property similar to the results of  Theorem  \ref{Ieq-33}, then one can prove that such algorithm converges linearly provided that the weighted matrix $G$ is positive definite.


\begin{thebibliography}{99}

\footnotesize
{
\bibitem{BaiLiXuZhang2017}
Bai, J.C.,    Li, J.C.,  Xu, F.M., Zhang, H.C.:
Generalized symmetric ADMM for separable convex optimization.
Comput. Optim. Appl.     70, 129-170 (2018)

\bibitem{Eckstein94}
Eckstein, J.:
Some saddle-function splitting methods for convex programming.
Optim. Methods Softw.  4,   75-83   (1994)


\bibitem{FPang2003}
 Facchinei, F.,   Pang, J.S.:
 Finite-Dimensional Variational Inequalities and Complementarity Problems. Springer-Verlag, Berlin  (2003)

\bibitem{Fazelun13}
Fazel, M.,  Pong, T.K.,  Sun, D.F.,   Tseng, P.:
Hankel matrix rank minimization with applications to system identification and realization.
SIAM J. Matrix Anal. Appl. 34, 946-977  (2013)

\bibitem{gaoma18}
Gao, B., Ma, F.:
Symmetric alternating direction method with
indefinite proximal regularization for linearly
constrained convex optimization.
J. Optim. Theory Appl.  176, 178-204 (2018)

\bibitem{GlowinskiMarrocco1975}
  Glowinski, R.:  Marrocco, A.:
Approximation par$\acute{e}$l$\acute{e}$ments finis d'rdre un  et r$\acute{e}$solution, par p$\acute{e}$nalisation-dualit$\acute{e}$  d'une classe de probl$\grave{e}$mes de Dirichlet non lin$\acute{e}$aires.
 Rev. Fr. Autom. Inform. Rech. Op$\acute{e}$r. Anal. Num$\acute{e}$r. 2,  41-76 (1975)


\bibitem{Heyuan2015}
 He, B.S.,   Yuan, X.M.:
 Block-wise alternating direction method of multipliers for multiple-block convex programming and beyond.
 SMAI J.  Comput. Math. 1,  145-174 (2015)



\bibitem{HeMaYuan2016}
He, B.S.,    Ma, F.,    Yuan, X.M.:
Convergence study on the symmetric version of ADMM with larger step sizes.
SIAM J. Imaging Sci.  9,   1467-1501  (2016)

\bibitem{HeY16}
He, B.S., Xu, H.K., Yuan, X.M.:
On the proximal Jacobian decomposition of ALM for multiple-block
separable convexminimization problems and its relationship to ADMM.
J. Sci. Comput. 66, 1204-1217 (2016)


 \bibitem{Robinson1981}
Robinson, S.M.:
Some continuity properties of polyhedral multifunctions.
Math. Program. Stud. 14, 2016-241 (1981)

 \bibitem{SUn18}
Sun, M., Sun, H.C.:
Improved proximal ADMM with partially parallel
splitting for multi-block separable convex programming.
J. Appl. Math. Comput.  58, 151-18 (2018)

 \bibitem{SS18}
 Sun, H.C.,   Tian, M.Y.,  Sun, M.:
The symmetric ADMM with indefiniteproximal regularization and its application.
J. Inequal.  Appl. 2017:172 (2017)

\bibitem{XuWu11}
 Xu, M.H.,    Wu, T..:
A class of linearized proximal alternating direction methods.
J. Optim. Theory Appl.  151,321-337  (2011)

}

\end{thebibliography}
\end{document}